\documentclass{article}

\usepackage{amsmath}
\usepackage{amsfonts}
\usepackage{amssymb}
\usepackage{tikz-cd}
\usepackage{graphicx}
\usepackage{hyperref}

\usepackage{fancyhdr}

\newtheorem{theorem}{Theorem}[section]
\newtheorem{lemma}[theorem]{Lemma}

\newenvironment{proof}[1][Proof]{\noindent\textbf{#1.} }{\ \rule{0.5em}{0.5em}}

\begin{document}

\title{The set of hyperbolic equilibria and of invertible zeros on the unit ball is computable}
\author{Daniel S. Gra\c{c}a\\Universidade do Algarve, C. Gambelas, 8005-139 Faro, Portugal\\\& Instituto de Telecomunica\c{c}\~{o}es, Lisbon, Portugal
	\and N. Zhong\\DMS, University of Cincinnati, Cincinnati, OH 45221-0025, U.S.A.}
\maketitle

\begin{abstract}
    In this note, we construct an algorithm that, on input of a description of a
    structurally stable planar dynamical flow $f$ defined on the closed unit disk,
    outputs the exact number of the (hyperbolic) equilibrium points and their
    locations with arbitrary accuracy. By arbitrary accuracy it is meant that any
    accuracy required by the input can be achieved. The algorithm can be further extended
    to a root-finding
    algorithm that computes
    the exact number of zeros as well the location of each zero of a
    continuously differentiable function $f$ defined on the closed unit ball of
    $\mathbb{R}^{d}$, provided that the Jacobian of $f$ is
    invertible at each zero of $f$; moreover, the computation is uniform in $f$.
\end{abstract}

\section{Introduction}

\thispagestyle{fancy}
Consider the autonomous system of ordinary differential equations
\begin{equation}
\label{ODE_Main}\dot{x}=f(x)
\end{equation}
for $f\in C^{1}(K)$, where $C^{1}(K)$ is the set of all continuously
differentiable functions in an open subset of $\mathbb{R}^{d}$
containing $K$ with values in $\mathbb{R}^{d}$, and
$K=\overline{B}(0, 1)$ is the closed unit ball centered at the
origin in $\mathbb{R}^{d}$. An equilibrium (or equilibrium point) of
(\ref{ODE_Main}) is a solution that does not change with time; that
is, the system (\ref{ODE_Main}) has an equilibrium solution
$x(t)=x_{0}$ if and only if $f(x_{0})=0$, where $x_0\in K$ and
$0\in\mathbb{R}^d$.

An equilibrium is the simplest possible solution of a dynamical
system. Nevertheless, it is fundamentally important because the
equilibria form a basis for analyzing more complicated behavior.
Yet, finding an equilibrium or, equivalently, solving the equation
$f(x)=0$ is easy
only in a few special cases and, in general, the equilibria of (\ref{ODE_Main}
) cannot be located exactly but only to be approximated by numerical
root-finding algorithms. There are numerous root-finding algorithms
such as Newton's method, Bisection method, Secant method, and
Inverse Interpolation method, to mention just a few. A numerical
algorithm is usually efficient when it is applied to some special
classes of elementary functions such as polynomials or augmented
with additional information such as good initial guesses or stronger
smoothness of the function.

There are also topological algorithms for computing zeros (or,
equivalently, fixed points) of a $C^k$ function defined on a compact
subset of $\mathbb{R}^d$ under certain conditions (see, for example,
\cite{Mil02}, \cite{Col08c}, \cite[Proposition 4.7]{BLMP19},
\cite[Section 5.1]{Neu19} and references therein). The topological
algorithms usually compute the set of all zeros as a \textit{subset}
of $\mathbb{R}^d$ and, for this reason, are often unable to exhibit
the position of each zero and to provide the exact number of the
zeros, if the set of zeros is finite. In other words, a topological
algorithm can provide arbitrarily good adumbrations of the set of
zeros but may not be able to do so for each zero individually. On
the other hand, from the perspective of dynamical systems, what
matters is the nature and the location of each individual
equilibrium; for example,  the Hartman-Grobman theorem, an important
theorem in dynamical systems, shows that the nonlinear vector field
$f(x)$ is conjugate to its linearization $Df(x_0)x$ in a
neighborhood of a hyperbolic equilibrium $x_0$, where the
neighborhood does not contain any other equilibrium point(see
e.g~\cite{Per01}).

In this note, we construct an algorithm that computes the exact
number and the locations of the zeros of a continuously
differentiable function $f$ defined on a compact subset of
$\mathbb{R}^{d}$, provided that the Jacobian of $f$ is invertible at
each zero of $f$ and that none of these zeros is on the boundary of the compact set;
moreover, the computation is uniform in $f$.

We first present our algorithm for computing the exact number and
the locations of the equilibrium points of a structurally stable
planar dynamical system, for that is our motivation. Intuitively, a
dynamical system is structurally stable if the qualitative behavior
of its trajectories persists under small perturbations. More
formally (see e.g.~\cite[Definition 1 in p.~317]{Per01}), the system
(\ref{ODE_Main}) is structurally stable if there is some
$\varepsilon>0$ such that for any $g\in C^1(K)$, if  $\| f - g\|_{1,\infty}
<\varepsilon$ on $K$, then there exists a homeomorphism $h: K \to K$
such that $h$ maps every trajectory of (\ref{ODE_Main}) onto a
trajectory of $y^{\prime}=g(y)$ preserving time orientation, where
$\| f\|_{1,\infty} = \max_{x\in K}\|f(x)\| + \max_{x\in K}\|Df(x)\|$.

If the system (\ref{ODE_Main}) is structurally stable on $K$, then
the number of its equilibria won't change when $f$ is slightly
perturbed.  Moreover, the location of every equilibrium depends
continuously on the perturbation (see e.g.~\cite[Theorems 1 and 2 in
p.~321]{Per01}). The equilibrium points which resist the small
perturbations are called hyperbolic equilibria. More precisely, an
equilibrium $x_0$ of (\ref{ODE_Main}) is said to be hyperbolic if
all eigenvalues of the Jacobian matrix of $f$ at $x_0$ have non-zero
real parts. This implies that the Jacobian of $f$ is invertible at
each hyperbolic equilibrium. Hyperbolic equilibria are as robust as
expected: small perturbations on $f$ do not alter the topological
character of the phase portrait near a hyperbolic equilibrium, but
only distort the trajectories near it by a small amount according to
the Hartman-Grobman theorem. Furthermore, if (\ref{ODE_Main}) is a
structurally stable planar dynamical system defined on $K$ (with
$d=2$), then the system has only finitely many equilibrium points
and all of them are hyperbolic according to the celebrated Peixoto
theorem \cite{Pei59}. The Peixoto theorem also shows that being
structurally stable is a generic property in the sense that the set
of all structurally stable planar vector fields is an open set of
$C^{1}(K)$ and it is dense in $C^{1}(K)$. Peixoto's theorem however
does not contain \textit{quantitative} information such as the exact
number of the equilibrium points or their locations inside the open
unit disk. Our algorithm, on the other hand, will supply such
quantitative information. The following theorem is the main result
of the note; its proof is given in Section \ref{Sec:MainResult1}.

Let $K_2$ be the closed unit ball when $d=2$,
i.e.~$K_2=\overline{B}(0, 1)\subseteq\mathbb{R}^2$. Take
$\mathcal{S}(K_2) = \{ f\in C^1(K_2) : \mbox{the system defined by
(\ref{ODE_Main}) is structurally stable}\}$ and let $Zero(f)$ be the
set of zeros of $f$ in $K_2$; i.e., the set of the equilibrium
points of (\ref{ODE_Main}). We recall that a $C^{1}$-name of a
real-valued function $f$ defined on a compact
$\Omega\subseteq\mathbb{R}$ is a sequence $\{ P_{k}\}$ of
polynomials with rational coefficients satisfying $\| f -
P_{k}\|_{1,\infty} <2^{-k}$ on $\Omega$; $f$ is computable if there
is a Turing machine that outputs the coefficients of $P_k$ on input
$k$. A name for a vector-valued function defined on
$\Omega\subseteq\mathbb{R}^d$ consists of the names of its component
(real-valued) functions. A vector-valued function is computable if
each of its components is computable.
(We use $K$ to denote the closed unit ball in
$\mathbb{R}^d$ for different dimension $d$; the dimension may not be
specified if the context is clear.)

\begin{theorem}
[Main result]\label{Th:main} The operator that assigns to each $f\in
\mathcal{S}(K_2)$ the set $Zero(f)$ and the exact size of $Zero(f)$
is computable. More precisely, there is an algorithm $\mathcal{E}$,
when given any $C^{1}$-name of $f\in\mathcal{S}(K_2)$ and
$n_{0}\in\mathbb{N}$ as input, $\mathcal{E}$ produces an integer
$n\geq n_{0}$, a nonnegative integer $\#_0(f)$, and a list $C$ of
finitely many squares with rational vertices as output, such that

\begin{enumerate}
\item $\#_0(f)$ is the exact number of the equilibrium points of
(\ref{ODE_Main});

\item each square in $C$ contains exactly one equilibrium point and has side of length $1/n$.
Furthermore, any zero of $f$ (in $K_2$) is contained in a square of $C$.
In particular, this implies that $d_{H}(Zero(f),\cup C)\leq1/n$,
where $d_{H}(Zero(f),\cup C)$ is the Hausdorff distance between $Zero(f)$
and $\cup C$.
\end{enumerate}
\end{theorem}

We note that each square in $C$ can be viewed as a pixel. This
result is restricted to planar systems, where the dimension is $d = 2$,
which is the usual case when considering structurally stable systems
following Peixoto's theorem. Nevertheless, Theorem \ref{Th:main} can
be generalized to a zero-finding algorithm for functions in
$\mathcal{Z}(K)$, which works in all dimensions $d\geq 1$ (i.e.~for
all $C^1$ functions defined on $K=\overline{B}(0,
1)\subseteq\mathbb{R}^d$, $d\geq 1$) and is more general in the
sense that it can also compute non-hyperbolic zeros $x$ of $f$ for
which $Df(x)$ might have eigenvalues with zero real part, as long as
$\det Df(x)\neq 0$ (hyperbolic points satisfy this condition), where
\[
\mathcal{Z}(K)=\{ f\in C^{1}(K): \, \mbox{$\det Df(\alpha)\neq 0$
    and $\alpha\not\in \partial K$ whenever $\alpha$ is a zero of $f$}\},
\]
$\partial K$ denotes the boundary of $K$, and
$K\subseteq\mathbb{R}^d$ is the closed unit ball, $d\geq 1$. The
formal description of the zero-finding algorithm and its proof are
presented 
in Section \ref{Sec:MainResult2}. We note that if
(\ref{ODE_Main}) is structurally stable, then there is no
equilibrium on $\partial K$.

As we mentioned before, the continuity of $f$ alone is not strong
enough for being able to compute the zeros of f. Some additional
conditions are needed according to the following well-known facts in
computable analysis (see, for example, \cite{BHW08}): (1) there is a
computable function $f: [0, 1] \to\mathbb{R}$ having infinitely many
zeros but none is computable; (2) the zero set of a function $f$ is
not uniformly computable in $f$; i.e. there is no algorithm that
outputs an approximation of the zero set of $f$ with error less than
$1/n$ (or $2^{-n}$) in some measurements when taking as input a
description of $f$ and a natural number $n$; and (3) given that $f$
has only finitely many zeros, the number of zeros of $f$  is not
uniformly computable in $f$. In fact, (2) and (3) are true even for
families of elementary functions.

\section{Definitions}\label{Sec:Definitions}

Let $\| \cdot\|_{2}$ and $\| \cdot\|_{\infty}$ denote the Euclidean
norm and the maximum norm of $\mathbb{R}^{d}$, respectively; let $K$
be the closure of the open unit ball $B(0, 1)$ of $\mathbb{R}^{d}$
in $\| \cdot\|_{\infty}$ norm; let $C^{1}(K)$ be the set of
continuously differentiable functions defined on some open subset of
$\mathbb{R}^{d}$ containing $K$ with values in $\mathbb{R}^{d}$; and
let $Df(x)$ denote the (Fr\'{e}chet) derivative of $f$ at $x$ for
$x\in K$, which is a  linear transformation from $\mathbb{R}^{d}$ to
$\mathbb{R}^{d}$. Since $f$ is $C^{1}$, it follows that, for each
$x\in K$, the linear transformation $Df(x)$ is continuous (thus
bounded) and is the same as the linear transformation induced by the
$d\times d$ Jacobian matrix of the partial derivatives of $f$ at
$x$; moreover, $Df(x)$ is continuous in $x$.

Let $A=(a_{ij})$ be a $d\times d$ square matrix (of real entries). Then the
following matrix norms are all equivalent: $\| A\|_{2}=\sup_{y\in
\mathbb{R}^{d}, \| y\|_{2}=1}\| Ay\|_{2}=\sqrt{\rho(AA^{T})}$, $\|
A\|_{\infty} = \sup_{y\in\mathbb{R}^{d}, \| y\|_{\infty}=1} \| Ay\|_{\infty} =
\max_{i}\sum_{j=1}^{d}|a_{ij}|$, and $\| A\|_{HS} = \sqrt{\sum_{i, j=1}
^{d}a_{ij}^{2}}$, where $\rho(A)=\max_{1\leq j\leq d}|\lambda_{j}|$,
$\lambda_{j}$ is an eigenvalue of $A$. The first two norms are
called operator norms while the third is called the Hilbert-Schmidt
norm. \\

\noindent\textbf{Convention} Since the norms $\| \cdot\|_{2}$ and $\|
\cdot\|_{\infty}$ imposed on $\mathbb{R}^{d}$ are equivalent, and the three
norms defined on the square $d\times d$ matrices $A$ are also equivalent, in
what follows we use $\| \cdot\|$ to denote either of the five norms when there
is no confusion in the context.

Next we recall several notions from computable analysis. For more
details and rigorous definitions the reader is referred to
\cite{Wei00} and references therein. Let $X$ be a set and let
$\delta:\Sigma^\omega\to X$ be a surjective map, where $\Sigma$ is a
finite set (the alphabet) and $\Sigma^\omega=\{w|
w:\mathbb{N}\to\Sigma\}$ is the set of all one-way infinite
sequences over $\Sigma$. In this case, $\delta$ is called a
representation of the set $X$. If $\delta(w)=x$, then $w$ is called
a ($\delta$-)name of $x$; an element $x \in X$ is called
($\delta$-)computable if it has a computable name.

\section{Proof of the main result}\label{Sec:MainResult1}

Before giving a precise description of the algorithm $\mathcal{E}$
and proving Theorem \ref{Th:main}, we present a number of lemmas and
auxiliary results. In what follows we assume that
$f\in\mathcal{S}(K_2)$.

\begin{lemma}\label{lem:def_M}
    Let $M=\max_{x\in K_2}\{\Vert
    Df(x)\Vert, |\det(Df(x))|,1\}$, where $\| \cdot\|$ is an operator
    norm. Then $M$ is computable from $f$.
\end{lemma}

\begin{proof}
Let $M=\max_{x\in K_2}\{\Vert
Df(x)\Vert, |\det(Df(x))|,1\}$, where $\| \cdot\|$ is an operator
norm. Then it follows from the discussion in the first paragraph of
the previous section that $Df(x)$ is a bounded linear
operator for each $x\in K_2$ and is continuous in $x$. Moreover, since a $C^{1}
$-name of $f$ is given as a part of the input to the algorithm to be
constructed, the Jacobian matrix of the partial derivatives of $f$
is computable uniformly in $x$ on the compact set $K_2$; thus the
map $\Vert Df\Vert:K_2\rightarrow\mathbb{R}^{+}$, $x\rightarrow\Vert
Df(x)\Vert$, is computable from $x$ and $f$ and, as a consequence,
$M$ is computable from $f$ (and $K_2$).
\end{proof}

\begin{lemma}\label{lem:upper_bound}
Let $\mathbb{D}$ be a small disk of $\mathbb{R}^2$ containing the
origin such that $x+h$ is in the domain of $f$ for every $x\in K_2$
and $h\in \mathbb{D}$. Then there is a decreasing sequence
$\{r_{m}\}_{m\geq1}$, computable from $f$, of positive numbers
$r_{m} \leq \min\{1, \mbox{the radius of } \mathbb{D}\}$ such that
whenever $\Vert h \Vert\leq r_m$, one has
\begin{equation}
\Vert f(x+h)-f(x)-(Df(x))(h)\Vert
\leq2^{-m-1}\Vert h\Vert
\label{e:continuity-of-Df}
\end{equation}

\end{lemma}

\begin{proof}
Let $g(t)=f(x+th)-f(x)-tDf(x)h$, where $t\in [0, 1]$.
 Then $g^{\prime}(t)=Df(x+th)h-Df(x)h=(Df(x+th)-Df(x))h$ and
 $g(0)=0$. Subsequently,
 \begin{align*}
 \left\Vert f(x+h)-f(x)-Df(x)h\right\Vert  & =\left\Vert g(1)\right\Vert \\
 & \leq\sup_{t\in\lbrack0,1]}\left\Vert g(0)+g^{\prime}(t)(1-0)\right\Vert \\
 & \leq\sup_{t\in\lbrack0,1]}\left\Vert Df(x+th)-Df(x)\right\Vert \cdot\left\Vert h\right\Vert.
 \end{align*}
(the first inequality is obtained by applying the mean value theorem
componentwise to $g$). Since $Df$ is computable and admits a
computable modulus of continuity, there is a sequence
$\{r_{m}\}_{m\geq1}$, computable from $f$, of positive numbers
$r_{m} \leq \min\{1, \mbox{the radius of } \mathbb{D}\}$ such that
$\left\Vert Df(x+th)-Df(x)\right\Vert\leq 2^{-m-1}$ whenever $\Vert
h \Vert\leq r_m$ and $t\in[0,1]$. Thus
$\sup_{t\in\lbrack0,1]}\left\Vert Df(x+th)-Df(x)\right\Vert\leq
2^{-m-1}$ whenever $\Vert h \Vert\leq r_m$, which implies that
for all $x\in K_2$ and all $\Vert h\Vert\leq r_{m}$, $\Vert h\Vert\neq0$.
Without loss of generality we may also assume that $r_{m}\geq r_{m+1}$. This concludes the proof of the lemma.
\end{proof}

\begin{lemma}\label{lem:lower_bound}
    Let $K_{Df}^{m}=\{x\in K_2:\,2^{-m}\leq\min\{\Vert
    Df(x)\Vert,|\det(Df(x))|\}\leq M\}$, where
    $m\in\mathbb{N}$ is arbitrary and $M$ is given in Lemma \ref{lem:def_M}. Let $\{r_{m}\}_{m\geq1}$ be the sequence defined in Lemma \ref{lem:upper_bound}. Then for each $x\in
    K_{Df}^{m}$ and every $\Vert h\Vert\leq r_{m}$, $\Vert h\Vert\neq0$,
    \begin{equation}
    \Vert f(x+h)-f(x)\Vert\geq2^{-m-1}\Vert h\Vert\label{e:lower-bound-on-f}
    \end{equation}
\end{lemma}

\begin{proof}
The estimate (\ref{e:lower-bound-on-f}) is obtained from the following calculation: since $x\in
K_{Df}^{m}$, it follows that $\Vert Df(x)\Vert\geq2^{-m}$, which implies that
\begin{align*}
2^{-m}\Vert h\Vert &  \leq\left\Vert (Df(x))(h)\right\Vert \\
&  =\left\Vert f(x+h)-f(x)-(Df(x))(h)+f(x)-f(x+h)\right\Vert \\
&  \leq\left\Vert f(x+h)-f(x)-(Df(x))(h)\right\Vert +\left\Vert
f(x)-f(x+h)\right\Vert \\
&  \leq2^{-m-1}\left\Vert h\right\Vert +\left\Vert
f(x)-f(x+h)\right\Vert
\end{align*}
Consequently
\begin{align*}
\left\Vert f(x)-f(x+h)\right\Vert  &  \geq2^{-m}\left\Vert h\right\Vert
-2^{-m-1}\left\Vert h\right\Vert \\
&  =2^{-m-1}\left\Vert h\right\Vert
\end{align*}
\end{proof}
We note that if $\det Df(x)\neq0$ for some $x\in K_2$, then $\Vert
Df(x)\Vert_{HS}\neq0$ which would then imply that $\Vert Df(x)\Vert_{\infty
}\neq0$. Therefore, if $Df(x)$ is invertible at some $x\in K_2$, then $\det
Df(x)\neq0$ and $\Vert Df(x)\Vert\neq0$.

The lemma below is also needed for the construction of the algorithm $\mathcal{E}$.
Its proof can be found in \cite{Rud06}.

\begin{lemma}
\label{Lemma:intermediate-value} Let $g:A\rightarrow\mathbb{R}^{d}$ be a
$C^{1}$ function defined on an open subset of $\mathbb{R}^{d}$ containing the
closed ball $A=\overline{B}(x_{0},r)$ with $r>0$; let $M=\max_{x\in
A}\left\Vert Dg(x)\right\Vert _{\infty} $. Then for all $x,y\in A$ we have
$\left\Vert g(y)-g(x)\right\Vert _{\infty} \leq M\left\Vert y-x\right\Vert
_{\infty} $.
\end{lemma}

Some notations are called upon before presenting the algorithm $\mathcal{E}$. Let
$C_{n}$ be a set of small squares
of side length $1/n$, called $n$-squares which covers exactly
$K_2$. More precisely,
let $G_{n}=((\mathbb{Z}+1/2)/n)^{2}\cap K_2$, where $(\mathbb{Z}
+1/2)/n=\{(z+1/2)/n\in\mathbb{R}:z\in\mathbb{Z}\}=\{\ldots,-\frac{3}
{2n},-\frac{1}{2n},\frac{1}{2n},\frac{3}{2n},\ldots\}$; let $C_{n}
=\{\overline{B}(x,1/(2n)):x\in G_{n})\}$, where
$\overline{B}(x,1/(2n))=\{y\in \mathbb{R}^{d}:\Vert
y-x\Vert_{\infty}\leq1/(2n)\}$.

We now define the algorithm $\mathcal{E}$. We will explain why the
algorithm works after presentation of the algorithm. Note that when
we can compute a quantity such as $d(0,f(s))$ in Step (3-1), we can
compute it with accuracy $2^{-(n+1)}$, obtaining an approximation
$r$. Therefore, if $r \leq 3/2\cdot 2^{-n}$, then we conclude that
$d(0,f(s))\leq 2^{-(n-1)}$. Otherwise, we have for sure
that $d(0,f(s))>2^{-n}$. We write this procedure as ``decide whether
$d(0,f(s))>2^{-n}$ or $d(0,f(s))\leq2^{-(n-1)}$''.

On input of (a
$C^{1}$-name of) $f$ and $n_{0}$ ($n_{0}\in\mathbb{N}$ defines the
accuracy $1/n_{0}$ of the result), the algorithm $\mathcal{E}$ runs
as follows:

\begin{enumerate}
\item Set $n=3n_{0}.$

\item For each $n$-square $s\in C_{n}$, set $\operatorname{result}(s)=$
Undefined and $\operatorname{counter}(s)=0$.

\item For each $n$-square $s\in C_{n}$, do the following:

\begin{itemize}
\item[(3-1)] $\mathcal{E}$ computes $d(0,f(s))$ until it decides whether
$d(0,f(s))>2^{-n}$ or $d(0,f(s))\leq2^{-(n-1)}$. If $d(0,f(s))>2^{-n}$, then
$\mathcal{E}$ sets $\operatorname{result}(s)=$ False and go to step (3-3).
Otherwise, compute $d(s,\partial K_2)$, increment $n$, and go to step 2 if
$d(s,\partial K_2)<5/n$.

\item[(3-2)] If $d(s,\partial K_2)\geq\frac{4}{n}$, then let $\mathcal{N}
(s)=\{x\in K_2:d(x,s)\leq1/n\}\subseteq K_2$ i .e. $\mathcal{N}(s)$ consists of
$s$ and the $n$-squares adjacent to $s$ and let $\mathcal{M}(s)=\{x\in
K_2:d(x,s)\leq3/n\}\subseteq K_2$. Clearly $s\subseteq\mathcal{N}(s)\subseteq
\mathcal{M}(s)\subseteq K_2$.$\ \mathcal{E}$ now computes $\min_{x\in
\mathcal{M}(s)}\{\Vert Df(x)\Vert,\left\vert \det(Df(x))\right\vert \}$ until
it decides whether $\min_{x\in\mathcal{M}(s)}\{\Vert Df(x)\Vert,\left\vert
\det(Df(x))\right\vert \}\leq2^{-n+1}$ or $\min_{x\in\mathcal{M}(s)}\{\Vert
Df(x)\Vert,\left\vert \det(Df(x))\right\vert \}>2^{-n}$. In the case the condition $\min
_{x\in\mathcal{M}(s)}\{\Vert Df(x)\Vert,\left\vert \det(Df(x))\right\vert
\}\leq2^{-n+1}$ holds, then increment $n$ and go to step 2.

\item[(3-3)] Repeat step (3-1) with a new $n$-square from $C_{n}$ or proceed
to step 4 if no $n$-square is left in $C_{n}$.
\end{itemize}

\item For each $n$-square $s\in C_{n}$, do the following:

\begin{itemize}
\item[(4-1)] If $\operatorname{result}(s)=$ False, then go to step (4-3).

\item[(4-2)] Pick $\tilde{n}\geq n$ such that for every $\tilde{n}$-square
$s_{j}\subseteq s$, it holds true that $\mathcal{M}(s_{j})\subset B(x_{s_{j}
},r_{n})$, where $x_{s_{j}}$ is the center of $s_{j}$ and
$B(x_{s_{j}},r_{n})$ is the open ball centered at $x_{s_{j}}$ with
radius $r_{n}$, where $r_n$ is given by Lemma \ref{lem:upper_bound}. Note that it follows from Lemma \ref{lem:lower_bound} and
(\ref{e:lower-bound-on-f}) that $f$ is injective on
$\mathcal{M}(s_{j})$. Let $s_{j}$, $1\leq j\leq J(s)$, be
$\tilde{n}$-squares such that $s=\bigcup_{j=1}^{J(s)}s_{j}$ and any
two distinct squares are either disjoint or intersect only in their
boundaries. Note that this condition holds only when $\tilde{n}=jn$
for some $j\in\mathbb{N}$. We assume without loss of generality that
this requirement is satisfied. Set $l=\tilde{n}$. For each $s_{j}$,
$1\leq j\leq J(s)$, do the following:

\begin{itemize}
\item[(4-2-1)] $\mathcal{E}$ computes $d(0,f(s_{j}))$. Set
$\operatorname{result}(j,s)=$ False if $d(0,f(s_{j}))>2^{-l}$; or go to
(4-2-2) if $d(0,f(s_{j}))\leq2^{-l+1}$.

\item[(4-2-2)] $\mathcal{E}$ computes rational points $x_{1},x_{2}
,\ldots,x_{e(j,l)}$ in the interior of $s_{j}$ such that $s_{j}\subseteq
\bigcup_{i=1}^{e(j,l)}B(x_{i},\frac{2^{-l-1}}{4M})$; afterwards, $\mathcal{E}$
computes the numbers $\theta_{i}=d(x_{i},\partial\mathcal{N}(s_{j}))$, $1\leq
i\leq e(j,l)$. We note that $\theta_{i}\leq r_{n}$ because $\mathcal{N}(s_j)\subset \mathcal{M}
(s_{j})\subset B(x_{s_{j}},r_{n})$ and $\theta_{i}\geq1/\tilde{n}$
because $x_{i}\in s_{j}$ and
$d(s_{j},\partial\mathcal{N}(s_{j}))\geq1/\tilde{n}$. Let
$B_{i}=B(x_{i},\theta_{i})$, $1\leq i\leq e(j,l)$. Then
$B_{i}\subset \mathcal{N}(s_{j})$. Next, $\mathcal{E}$ computes
$d_{i}=\min_{x\in\partial B_{i}}\Vert f(x)-f(x_{i})\Vert$ for each
$1\leq i\leq e(j,l)$. Since $\theta_{i}\leq r_{n}$ and
$\min_{x\in\mathcal{M}(s)}\{\Vert Df(x)\Vert,|\det
Df(x)|\}\geq2^{-n}$, it follows from (\ref{e:lower-bound-on-f}) that
$\Vert f(x)-f(x_{i})\Vert\geq2^{-n-1}\theta_{i}$ for every
$x\in\partial B_{i}$; thus
$d_{i}\geq2^{-n-1}\theta_{i}\geq2^{-n-1}/\tilde{n}$ for all $1\leq
i\leq e(j,l)$.

\item[(4-2-3)] Let $D_{i}=B(f(x_{i}),d_{i}/2)$.

\begin{enumerate}
\item If $0\in\bigcup_{i=1}^{e(j,l)}D_{i}$ and $s_{j}$ is adjacent to
an $\tilde{n}$-square $s_{k}^{\ast}$ (note that this
$\tilde{n}$-square may belong to an $n$-square $s^{\prime}\neq s$ if
$s_{j}$ is adjacent to the boundary of $s$) for which
$\operatorname{result}(k,s^{\ast})$ is defined and equal to True,
then $\mathcal{E}$ sets $\operatorname{result}(j,s)=$ False and
moves to the next $\tilde{n}$-square $s_{j+1}$.

\item If $0\in\bigcup_{i=1}^{e(j,l)}D_{i}$ and $s_{j}$ is not adjacent to
an $\tilde{n}$-square $s_{k}^{\ast}$ for which
$\operatorname{result}(k,s^{\ast })$ is defined and equal to True,
then $\mathcal{E}$ sets $\operatorname{result}(j,s)=$ True and moves
to the next $\tilde{n}$-square $s_{j+1}$.

\item If $0\not\in\bigcup_{i=1}^{e(j,l)}D_{i}$, then $\mathcal{E}$ increments $l$ by 1
and returns to (4-2-1).
\end{enumerate}

\item[(4-2-4)] Repeat step (4-2-1) with the next $\tilde{n}$-square $s_{j+1}$ or
proceed to step (4-2-5) if no more $\tilde{n}$-square is left in
$\{s_{1},s_{2},\ldots,s_{J(s)}\}$.

\item[(4-2-5)] Set $\operatorname{result}(s)=$ False if $\operatorname{result}
(j,s)=$ False for all $1\leq j\leq J(s)$. Otherwise, set
$\operatorname{result}(s)=$ True; set $\widetilde{\mathcal{N}}(s)=\bigcup
\{\mathcal{N}(s_{j})\,:\,\mbox{$\operatorname{result}(j, s) =$ True}\}$; and
set $\operatorname{counter}(s)=$ the cardinality of the set
$\{j:\,\mbox{$\operatorname{result}(j, s) =$ True},1\leq j\leq J(s)\}$.
\end{itemize}

\item[(4-3)] Repeat step (4-1) with a new $n$-square from $C_{n}$ or proceed
to step 5 if there is no more $n$-square left in $C_{n}$.
\end{itemize}

\item Let $C=\varnothing$.

\item For all $s\in C_{n}$ do:

\begin{enumerate}
\item If $\operatorname{result}(s)=$ True, then $C=C\cup\widetilde
{\mathcal{N}}(s)$.

\item Set $\#_0(f)=\sum_{s\in C_{n}}\operatorname{counter}(s)$.
\end{enumerate}

\item Output $C$ and $\#_0(f)$.
\end{enumerate}

\noindent\textbf{Proof of Theorem \ref{Th:main}.} We now show that the algorithm $\mathcal{E}$ works. In step 1, the
accuracy is increased because $\mathcal{E}$ may not be able to tell
whether a zero is in an $n$-square $s$ but, nevertheless, it is
capable of determining if a zero is in a $(1/n)$-neighborhood of
$s$. Therefore, in order to get the requested accuracy bound of
$1/n_{0}$, we will have to use smaller squares whose side lengths
are at most diameter $1/(3n_{0})$. This is why $\mathcal{E}$ starts
at $n=3n_{0}$.

After step 1, $\mathcal{E}$ studies whether or not $s$ contains a
zero of $f$ for every $s\in C_{n}$. The result of this investigation
is saved into $\operatorname{result}(s)$. This is why
$\operatorname{result}(s)$ is initially set to be undefined in step
2.

In step (3-1) $\mathcal{E}$ performs a first test (with accuracy
$2^{-n}$) to check whether $f(s)$ contains a zero. Clearly, a False
output at step (3-1) indicates that $s$ contains no zeros of $f$. If
this is the case $\mathcal{E}$
proceeds to another $n$-square. If we instead obtain $d(0,f(s))\leq2^{-(n-1)}
$, then $s$ may or may not have a zero, and so further
investigations are needed. Meanwhile $\mathcal{E}$ runs a test to
determine whether $d(s,\partial K_2)<5/n$ or $d(s,\partial K_2)\geq4/n$
for the purpose of getting rid of those $n$ squares which are too
close to the boundary $\partial K_2$ of $K_2$ because the further
investigations are to be performed on some $(3/n)$-neighborhood of
$s$. In the case that $d(s,\partial K_2)<5/n$, $\mathcal{E}$
increments $n$ and then repeats step (3-1). It follows from Lemma
\ref{Lem:FiniteNumberZeros} that $f$ has only finitely many zeros
and none is on $\partial K_2$; therefore, when $n$ is large enough,
$\mathcal{E}$ will output $d(0,f(s))\geq2^{-n+1}$ whenever
$d(s,\partial K_2)<5/n$ for every $s\in C_{n}$. In other words,
$\mathcal{E}$ will not go into a loop; it will either output
result(s)=False or it will move on to the next step for sufficiently
large $n$; for example, for $n$ that satisfies the condition:
$2^{-n+1}<\min\{\Vert f(x)\Vert :\,\mbox{$x\in K_2$ and $d(x,
\partial K_2)\leq \sigma/2$}\}$, where $\sigma=\min\{d(x,\partial K_2):\,f(x)=0\}$.

In step (3-2) $\mathcal{E}$ tests whether the jacobian $Df$ is
invertible on $\mathcal{M}(s)$ for those $s\in C_{n}$ satisfying
$d(s,\partial K_2)\geq4/n$; i.e., those $n$ squares whose final status
result(s) are in need to be updated
from Undefined to True or False. In case that $\min_{x\in\mathcal{M}(s)}
|\det(Df(x))|>2^{-n}$, the inverse function theorem can be applied to
determine whether or not $s$ contains a zero of $f$. And so we wish to put
each $n$-square $s$ satisfying $d(s,\partial K_2)\geq4/n$ into one of the two
groups: either $d(0,f(s))>2^{-n}$ or $\min_{x\in\mathcal{M}(s)}\min\{\Vert
Df(x)\Vert,|\det(Df(x))|\} >2^{-n}$. This is achievable because, once again
due to Lemma \ref{Lem:FiniteNumberZeros}, $f$ has only finitely many zeros,
say $\xi_{1},\ldots,\xi_{k}$, with the property that $|\det(Df(\xi_{i}))|>0$
and $\Vert Df(\xi_{i})\Vert>0$ as well (recall that $\Vert Df(x)\Vert$ is
equivalent to $\Vert Df(x)\Vert_{HS}$). By continuity of $Df$, there exist
$\rho_{1},\ldots,\rho_{k}>0$ such that $\min\{\Vert Df(y)\Vert,\left\vert
\det(Df(y))\right\vert \}>0$ for all $y$ in the closed ball $\overline{B}
(\xi_{i},\rho_{i})$, $1\leq i\leq k$. Then $\min\{\Vert Df(y)\Vert,\left\vert
\det(Df(y))\right\vert \}>0$ for all $y\in\overline{B}(\xi_{i},\rho)$, $1\leq
i\leq k$, where $0<\rho=\min_{1\leq i\leq k}\rho_{i}$. Thus when $n$ is large
enough meeting the conditions $4/n\leq\rho$ and
\[\min_{y\in\cup_{i=1}
^{k}\overline{B}(\xi_{i},\rho)}\min\{\Vert Df(y)\Vert,\left\vert
\det(Df(y))\right\vert \}\geq2^{-n+1},
\]
 the test of step (3-2) is guaranteed
to succeed for all squares $s$ which reach step (3-2).

In step 4 we fix the first $n$ which successfully led us to this step. We
observe that if $\operatorname{result}(j,s)=$ True, then $\mathcal{N}(s_{j})$
contains one and only one equilibrium because for every $1\leq i\leq e(j,s)$,$B_{i}
\subset\mathcal{N}(s_{j})\subset \mathcal{M}(s_j)$,  $f:B_{i}\bigcap
f^{-1}(D_{i})\rightarrow D_{i}$ is a homeomorphism (see, for
example, \cite{Spi65}), and $f$ is injective on $\mathcal{M}(s_j)$.
However there remains a potential problem of a zero of $f$ being
counted multiple times. This may happen when $\mathcal{N}(s_j)$
contains a zero of $f$ and $s_j$ is adjacent to a portion of a
common side shared by $s$ and another $n$-square $\tilde{s}$. If
$\mathcal{E}$ acts on $\tilde{s}$ before it picks up $s$, then the
zero of $f$ contained in $\mathcal{N}(s_j)$ may already be detected
and counted by $\mathcal{E}$ at the time when $\mathcal{E}$ was
working on $\tilde{s}$. The potential multiple-counting may also
happen when the interior of $\mathcal{N}(s_j)$ intersects the
interior of $\mathcal{N}(s_{j^{\prime}})$ for some
$\tilde{n}$-square $s_{j^{\prime}}$ (contained in the $n$-square
$s$) that is adjacent to $s_j$. If the intersection contains an
equilibrium and if $\mathcal{E}$ acts on $s_{j^{\prime}}$ prior to
picking up $s_j$, then the equilibrium would have been counted by
$\mathcal{E}$ at the time while working with $s_{j^{\prime}}$, if
not earlier. The step (4-2-3) is a preventive mechanism designed to
ensure that an equilibrium is counted exactly once. It is also true
that, for each $n$-square $s$ satisfying
$\min_{x\in\mathcal{M}(s)}\Vert Df(x)\Vert >2^{-n}$, $\mathcal{E}$
will halt and produce either result(s) = False or result(s) = True.
To see this let us fix an $s_{j}$ with $1\leq j\leq J(s)$. If
$s_{j}$ does not contain an equilibrium, then the inequality
$d(0,f(s_{j}))>2^{-l}$ would appear for $l$ large enough. On the
other hand, assume that $s_{j}$
contains an equilibrium $x_{0}$. Let us pick some $l$ such that $2^{-l-1}
<2^{-n-2}/\tilde{n}$. Then there is some $x_{i}$ in the interior of
$s_{j}$, $1\leq i\leq e(j,l)$, such that $\Vert
x_{i}-x_{0}\Vert\leq\frac{2^{-l-1}}{4M}$. It now follows from Lemma
\ref{Lemma:intermediate-value} that $\Vert
f(x_{0})-f(x_{i})\Vert\leq4M\Vert x_{0}-x_{i}\Vert\leq2^{-l-1}$; in
other words, $0\in B(f(x_{i}),2^{-l-1})\subseteq
B(f(x_{i}),d_{i}/2)=D_{i}$ (recall that
$d_{i}\geq2^{-n-1}\theta_{i}\geq2^{-n-1}/\tilde{n}$ for all $1\leq
i\leq e(j,l)$).

\section{Computing invertible zeros}\label{Sec:MainResult2}

It is shown in \cite[Theorem 6.3.2]{Wei00} that the multi-valued function $f\in
C[0, 1] \mapsto \{ (f, x): \, f(x)=0\}$ is not continuous and thus
not computable. This raises the question: are there
topological/regularity conditions which we can impose on the family
of functions to ensure the uniform computability of the zero sets
and the cardinalities of the zero sets. Here the answer is yes. The
problem is to find the right conditions.

In this section, we consider the compact set $K=\overline{B}(0,
1)\subseteq\mathbb{R}^d$, $d\geq 1$. The algorithm $\mathcal{E}$
suggests that the following conditions can be imposed on the family
of continuously differentiable functions which will ensure the
uniform computability of the zero sets and their sizes: let
\[
\mathcal{Z}(K)=\{ f\in C^{1}(K): \, \mbox{$\det Df(\alpha)\neq 0$ and $\alpha\not\in \partial K$
whenever $\alpha$ is a zero of $f$}\}
\]
where
$\partial K$ denotes the boundary of $K$. We note that if (\ref{ODE_Main}) is
structurally stable, then there is no equilibrium on $\partial K$. And if we
dismiss the condition that $\alpha\not \in \partial K$, the uniform
computability may no longer be guaranteed. For example, let $\mathcal{A}=\{
f_{a}(x) : \, a\in\mathbb{R}\}$, where $f_{a}(x)=x-a$, $x\in K=[0, 1]$. Then
$f_{a}\in C^{1}(K)$ and $f^{\prime}_{a}=1$ on $K$ for all $a$; but the
operator $\mathcal{A} \to\{ 0, 1\}$, $f_{a}\mapsto\#$ of zeros of $f_{a}$, is
not uniformly computable because it is not even continuous. The problem is
caused by two ``bad" functions $f_{0}$ and $f_{1}$ - the zero of $f_{0}$ and
the zero of $f_{1}$ lie on $\partial K$. If we get rid of these two functions,
then we can compute the number of zeros of $f_{a}$ uniformly on $\mathcal{A}
\setminus\{ f_{0}, f_{1}\}$.

The algorithm $\mathcal{E}$ described in the previous section can be
adapted in a straightforward manner to $K\subseteq\mathbb{R}^d$ for
$d\geq 1$ and be applied, together with the lemma below, to compute
the exact number and the locations of the zeros of
$f\in\mathcal{Z}(K)$, uniformly on $\mathcal{Z}(K)$, thus showing
Theorem \ref{Th:main-zero}.

\begin{lemma}
\label{Lem:FiniteNumberZeros}Let $f\in\mathcal{Z}(K)$. Then $f$ has at most a
finite number of zeros in $K$.
\end{lemma}

\begin{proof}
Assume otherwise that $f$ has infinitely many zeros $\alpha_{n}$,
$n\in\mathbb{N}$, in $K$. Then since $K$ is compact, it follows that
$\{\alpha_{n}\}$ has a convergent subsequence, say $\{\alpha_{n_{k}}\}$, that
converges to $\alpha\in K$, which results in the following limit
\[
0=\lim_{k\rightarrow\infty}f(\alpha_{n_{k}})=f(\lim_{k\rightarrow\infty}
\alpha_{n_{k}})=f(\alpha).
\]
Note that the limit can be taken into $f$ is due to the fact that $f$ is
continuous on a compact set. Thus $\alpha$ is a zero of $f$. Moreover, since
$f\in\mathcal{Z}(K)$, it follows that $Df(\alpha)$ is invertible, which is
equivalent to the condition that the Jacobian determinant of the matrix
$Df(\alpha)$ is nonzero. Therefore $f$ itself is invertible in a neighborhood
of $\alpha$. But this is a contradiction because in any neighborhood $B$ of
$\alpha$ there is an $\alpha_{n_{k}}\in B$ and thus $f(\alpha_{n_{k}
})=f(\alpha)=0$, which implies that $f$ cannot be injective in any
neighborhood of $\alpha$, no matter how small it is.
\end{proof}

The construction of the algorithm $\mathcal{E}^{\prime}$ is the same as that
of the algorithm $\mathcal{E}$.

\begin{theorem}[Computing the set of invertible zeros]\label{Th:main-zero}
    The operator that assigns to each $f\in\mathcal{Z}(K)$ its zero set and the
    cardinality of the zero set is computable. More precisely, there is an
    algorithm $\mathcal{E}^{\prime}$, when given any $C^{1}$-name of $f$ and
    $n_{0}\in\mathbb{N}$ as input, $\mathcal{E}^{\prime}$ produces an integer
    $n\geq n_{0}$, a nonnegative integer $\#_0(f)$, and a list
    $C$ of finitely many squares with rational vertices
    (or hyper-cubes if $d\geq3$) as output, such that

    \begin{enumerate}
    \item $\#_0(f)$ is the exact number of zeros of $f$;

    \item each square in $C$ has side length $1/n$ and contains exactly one zero
of $f$. Furthermore, $Zero(f)=\{x\in K:f(x)=0\}\subseteq \cup C$,
which implies that $d_{H}(Zero(f),\cup C)\leq1/n$, where
$d_{H}(Zero(f),\cup C)$ is the Hausdorff distance between $Zero(f)$
and $\cup C$.
\end{enumerate}
\end{theorem}

\textbf{Acknowledgments.} D. Gra\c{c}a was partially funded by FCT/MCTES
through national funds and co-funded EU funds under the project
UIDB/50008/2020. \includegraphics[width=3.5mm]{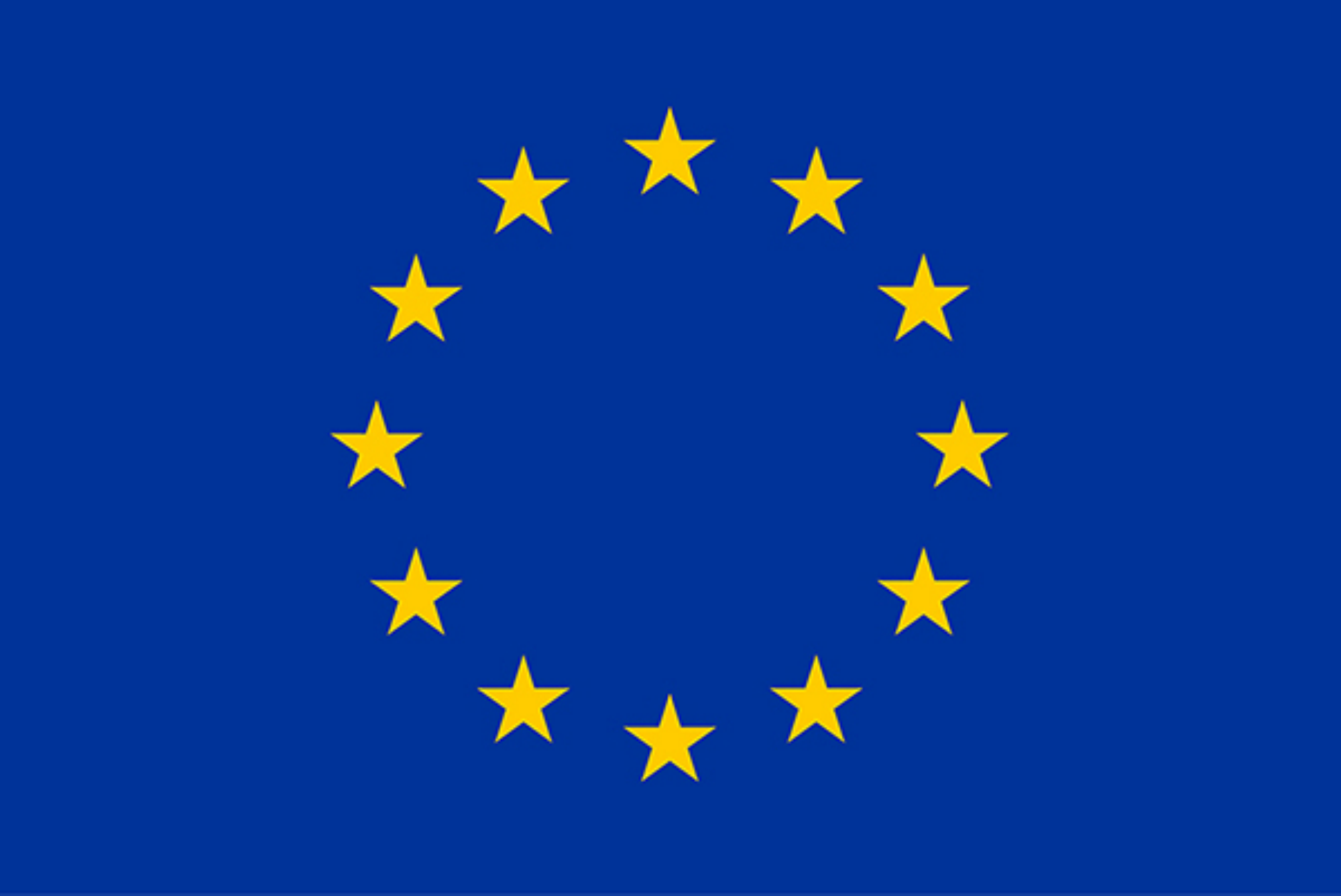} This project
has received
funding from the European Union's Horizon 2020 research and innovation
programme under the Marie Sk\l {}odowska-Curie grant agreement No 731143.

\bibliographystyle{amsplain}

\end{document}